\documentclass[multphys,vecphys]{svmult}

\usepackage{makeidx}         
\usepackage{graphicx}        
\usepackage{multicol}        
\usepackage[bottom]{footmisc}

\makeindex             

\begin{document}

\title*{Statistics of Extreme Spacings in Determinantal Random Point 
Processes}
\titlerunning{Statistics of Extreme Spacings} 
\author{Alexander Soshnikov \thanks{Dedicated to Yakov Sinai with great respect and admiration on the occasion of his 70th birthday}
}
\institute{University of California at Davis\\
Department of Mathematics\\
Davis, CA 95616, USA\\
\texttt{soshniko@math.ucdavis.edu}}
\maketitle

\section{Introduction}
\label{sec:1}

Determinantal (a.k.a. fermion) random point processes were introduced in probability theory by Macchi about thirty years ago (\cite{Ma1}, 
\cite{Ma2}, \cite{DVJ}). In the last ten years the subject has attracted a considerable attention due to its rich connections to 
Random Matrix Theory, Combinatorics, Representation Theory, Random Growth Models, Number Theory and several other areas of mathematics.
We refer the reader to the recent surveys (\cite{S1}, \cite{L}, \cite{HKPV}),
and research papers on the subject (\cite{BG}, \cite{DE}, \cite{GY}, \cite{Jo1}, \cite{Jo2}, \cite{LS}, 
\cite{Lyt}, \cite{PV}, \cite{ST1}, \cite{ST2},  \cite{SY}, \cite{S3}, \cite{S4}, \cite{Yoo}). 

In this paper we shall consider determinantal random point processes  on the real line with the translation-invariant correlation kernel.
In other words,  a one particle space $X$ is given as $X=\mathcal{R}^1, \ $ and the space of elementary outcomes $\Omega \ $  consists
of the countable, locally finite particle configurations on the real line
$$ \Omega = \{ \xi=(x_i)_{i=-\infty}^{+\infty} : \ \ \#( x_i \in [-N, N]) <+\infty, \ \ \forall N>0 \}, $$
where
$ x_i \in \mathcal{R}^1, \ i=0, \pm 1, \pm 2, \ldots, $ and $ \ \#( x_i \in [-N, N])=:\#([-N, N]) \ $ denotes the number of the 
particles in the 
interval $[-N, N]. \ $
Let us denote the set of the non-negative integers by $\mathcal{Z}_{+}^1= \{0,1,2,\ldots\}, \ $ and the set of the positive integers by 
$\mathcal{N}=\{1,2,\ldots\}. \ $
The $\sigma$-algebra $ \ \mathcal{F}$ of the measurable subsets of $\Omega$ is generated by the cylinder sets
$ \ C_{I_1, I_2, \ldots, I_k}^{n_1, n_2, \ldots n_k}= \{ \xi: \ \#(I_j)=n_j, \ j=1, \ldots, k\}, \ $ where $k$ is an arbitrary positive 
integer, $ \ k\in \mathcal{N}, \ $
$I_1, \ldots, I_k \ $ are arbitrary disjoint 
subintervals of the real line,  and $n_1, n_2, \ldots, n_k \in \mathcal{Z}^1_+.$

A probability measure $\mathcal{P}$ on the measurable space $(\Omega, \mathcal{F}) \ $ defines a random point process 
$(\Omega, \mathcal{F}, \mathcal{P}).\ $ A random point process is called determinantal if its $k$-point correlation functions have 
determinantal 
form
\begin{equation}
\label{detcorrfun}
\rho_k(x_1,\ldots, x_k)=\det(K(x_i,x_j))_{i,j=1, \ldots, k}, \ \ k=1,2,\ldots,
\end{equation}
where $K(x,y)$ is usually called the correlation kernel of the determinantal random point process.
We remind the reader that $k$-point correlation functions are defined in such a way that
\begin{equation}
\label{defcorr}
E \prod_{l=1}^k  \#(I_l) = \int_{I_1\times \ldots I_k} 
\rho_k(x_1, \ldots, x_k) \* dx_1\*\ldots dx_k
\end{equation}
for the arbitrary disjoint intervals $I_1, \ldots, I_k. \ $
Since the r.h.s. of (\ref{detcorrfun}) is non-negative, it follows that the correlation kernel $K(x,y)$ has non-negative minors.
If, in addition, the integral operator $K: L^2(\mathcal{R}^1) \to L^2(\mathcal{R}^1), \ \ (Kf)(x)= \int_{-\infty}^{+\infty} 
K(x,y)\*f(y)\* dy, $  
is Hermitian, one can conclude that $K$ is non-negative 
definite, i.e. $Spec(K) \in [0, +\infty).\ $ In the Hermitian case one can show that the necessary and sufficient condition on $K$ to define
a determinantal random point field (\ref{detcorrfun}) is 
\begin{equation}
\label{0K1}
0\leq K \leq 1, 
\end{equation}
in other words both $K$ and $1-K$ must be non-negative definite
operators (\cite{S1}, \cite{Ma1}). 

In this paper we consider the translation-invariant kernel 
\begin{equation}
\label{transl}
K(x,y)=g(y-x),\  \  \texttt{where}  \ \ g(x)=\int_{-\infty}^{+\infty} 
\exp(2\*\pi\* i \*x \*t) \* \phi(t)\* dt, 
\end{equation}
and $\phi(t) $ is an even real-valued integrable function
\begin{equation}
\label{fi}
\phi(t)=\phi(-t), \ \ \ \phi \in L^1(R^1).
\end{equation}
It follows from (\ref{0K1})
that 
\begin{equation}
\label{0fi1}
0 \leq \phi(t) \leq 1 \ \ (a.e.).
\end{equation}
In addition, we assume that the following technical conditions are satisfied
\begin{eqnarray}
\label{raz}
& & \int_{-\infty}^{+\infty} t^2\*\phi(t)\* dt < +\infty, \\
\label{dva}
& & |g(x)|\leq \frac{C}{1 + |x|^{\frac{1}{2} +\epsilon}}, \\
\label{tri}
& & | g'(x)|\leq \frac{C}{1 + |x|^{\frac{1}{2} + \epsilon}},
\end{eqnarray}
where $ \ C \ $ is a positive constant, and $ \ \epsilon \ $ is an arbitrary small positive constant.

Let $L$ be a large positive number. Consider a restriction of a configuration $\xi$ to the interval $[0,L].$ Let us denote the points of
$\xi \cap [0,L] \ $ by $x_1, x_2, \ldots, x_{N(L)}, \ $ where $N(L)$ is the cardinality of $\xi \cap [0,L].\ $
We assume that the points in $\xi \cap [0,L] \ $ are ordered:
$x_1 < x_2 < \ldots < x_{N(L)}. \ $  It is a well known (see e.g. \cite{S1}) that with probability 1 no two particles of a determinantal 
random point process coincide. We are interested to study the nearest spacings $ \theta_i= x_{i+1}-x_i, \ \ i=1, \ldots, N(L)-1,$
between the neighboring particles.  Functional Central Limit Theorem for the empirical distribution function of the 
nearest spacings of particles in $[0,L]$ (in the limit $L\to \infty$) was proven in \cite{S2} for $K(x,y)= \frac{\sin(\pi\*x)}{\pi\*x} $
(i.e. $\phi$ is the indicator of $[-\frac{1}{2},\frac{1}{2}]$), 
and for similar kernels arising in Random Matrix Theory. It was shown in \cite{S1} that 
the result could be extended to a quite general class of translation-invariant correlation kernels.

In this paper we study the smallest nearest spacings.
Our main result is the following

\begin{theorem}
Let $(\Omega, \mathcal{F}, \mathcal{P}) $ be a determinantal random point process on the real line with the translation-invariant correlation 
kernel $K(x,y)=g(y-x) $ satisfying
(\ref{transl})-(\ref{tri}).
Then the number of the nearest spacings less than $s/L^{1/3}$ in the interval $ [0,L] \ $ converges in distribution  
to the Poisson random variable with the  mean $\alpha\* s^3, $ in the limit $ L \to \infty, $  where
\begin{equation}
\label{nj}
\alpha=\frac{1}{3}\* g(0)\*g''(0)=
\frac{4 \* \pi^2}{3}\* \int_{-\infty}^{+\infty}  \phi(t) \* dt \* \int_{-\infty}^{+\infty} t^2\* \phi(t) \* dt.
\end{equation}
\end{theorem}

Let 
$ \ \eta(L)= L^{1/3} \* \min \{ \theta_i, \ \ i=1, \ldots, N(L)-1\}. \ $
In other words, $\eta(L)$ is the smallest nearest spacing in $[0,L],$ 
rescaled by $L^{1/3}.\ $
Theorem 1 immediately implies

\begin{theorem}
Let the conditions in Theorem 1 be satisfied.
Then
\begin{equation}
\label{limitminsp}
\lim_{L\to\infty} \Pr( \eta(L)>s)=\exp(-\alpha\*s^3).
\end{equation}
\end{theorem}

The method of the proof of Theorems 1 and 2 relies on the detailed analysis of $k$-point correlation and cluster functions
of the $s$-modified random point process, introduced in \cite{S1}, \cite{S2}. 

The rest of the paper is organized as follows. Point correlation and cluster  functions, and $s$-modified random point processes are 
discussed in Section 2. The proofs of Theorems 1 and 2  are given in Section 3.

We will use the notations $const,\ const_k, \ Const,  $ to denote various positive constants throughout this text. 
The values of these constants may be different in various 
parts of the paper. We shall use the notation $f=O(g) $ if the ratio $ f/g $ is bounded from above and below by some 
positive constants, and the notation $f=o(g) $ if the ratio $ f/g $ goes to zero.

\section{Correlation and Cluster Functions}
\label{sec:2}

We start this section by recalling the definition of a $k$-point cluster function
(sometimes also known as the Ursell factor).
For additional information we refer the reader to \cite{R}, \cite{LP}, \cite{CL}, \cite{S2}.

\noindent{\bf Definition}. {\it The $l$-point cluster function 
$r_l (x_1,\dots ,x_l ), \ l =1,2,\dots ,$ of a random point field
is defined in terms of the point correlation functions by the formula
\begin{equation}
\label{cluster}
r_l (x_1,\dots ,x_l )=\sum_G(-1)^{m-1}(m-1)!\*\prod^m_{j=1}
\rho_{\vert G_j\vert}(\bar x(G_j))
\end{equation}
where the sum is over all partitions $G$ of $[l ]=\{1,2,\dots ,l\}$
into subsets $G_1,\dots ,G_m,\ m=1,\dots ,l$, and $\bar x(G_j)=\{x_i:
i\in G_j\}, \ \vert G_j\vert =\#(G_j)$.}

The point correlation functions can be expressed in terms of the point cluster functions as

\begin{equation}
\label{corr-cherez-clust}
\rho_l (x_1,\dots ,x_l )=\sum_G\prod^m_{j=1}r_{\vert G_j\vert}
(\bar x(G_j)).
\end{equation}
The reader can observe that (\ref{cluster}) is the M\"obius inversion formula applied to (\ref{corr-cherez-clust}).
In particular, 
\begin{eqnarray}
\nonumber
& & \rho_1(x)= r_1(x),\\ 
\nonumber
& & \rho_2(x_1, x_2)= r_2(x_1, x_2) + r_1(x_1)\*r_1(x_2), \\ 
\nonumber
& & \rho_3(x_1, x_2, x_3)=
r_3(x_1, x_2, x_3) + r_2(x_1, x_2)\*r_1(x_3) +r_2(x_1, x_3)\*r_1(x_2) +\\
\nonumber
& &r_2(x_2, x_3)\*r_1(x_1) + r_1(x_1)\*r_1(x_2)\*r_1(x_3). 
\end{eqnarray}

It follows from (\ref{corr-cherez-clust}) and (\ref{detcorrfun}) that for  determinantal random point fields 
\begin{equation}
\label{davis}
r_l (x_1,\dots ,x_l )=(-1)^{l -1}\sum_{\texttt{cyclic }\sigma\in
S_l}K(x_1,x_{\sigma(1)})\*K(x_2,x_{\sigma(2)})\ldots  K(x_l ,x_{\sigma(l)}),
\end{equation}
where the sum in (\ref{davis}) is over all cyclic permutations.   In other words, for determinantal random  processes 
the difference between the formula (\ref{davis})
for the $l$-point cluster function and the formula 
\begin{equation}
\label{winters}
\rho_l (x_1,\dots ,x_l )=\sum_{\sigma\in S_l}(-1)^\sigma K(x_1,
x_{\sigma (1)})\cdot K(x_2,x_{\sigma (2)})\cdot\ldots\cdot K(x_l,x_{\sigma
(l)})
\end{equation}
for the $l$-point correlation function is that in (\ref{winters}) the summation is taken over all permutations in 
$ S_l, $ and in (\ref{davis}) the summation is over the cyclic permutations only.
One
can rewrite (\ref{davis}) as
\begin{equation}
\nonumber
r_l (x_1,\dots ,x_l )=(-1)^{l -1}\cdot\frac{1}{l}\sum_{\sigma
\in s_l}K(x_{\sigma (1)},x_{\sigma (2)})\*
K(x_{\sigma (2)},x_{\sigma (3)})\ldots K(x_{\sigma (l)},x_{\sigma
(1)}).
\end{equation}

It follows from (\ref{defcorr}) that the integral of
the $k$-point correlation function over the $k$-dimensional cube $[0,L]^k$ is equal to the 
$k$-th factorial moment of the counting random variable $\#(I), \  I=[0,L],\ $  namely
\begin{equation}
\label{tel}
E \* \#(I) \* (\#(I)-1) \ldots (\#(I)-k+1) = \int_{I^k} \rho_k(x_1, \ldots, x_k) \* dx_1\ldots dx_k.
\end{equation}

The integral of the $k$-point cluster function,
in turn, can be expressed as a linear combination of the cumulants of $\ \#([0,L]).\ $ Namely, let $C_k(L)$ denote the $k$-th cumulant of
$\#([0,L]), \ $ and $ \ V_k(L)= \int_{[0,L]^k} r_k(x_1,x_2,\ldots x_k) \* dx_1 \* dx_2 \ldots dx_k. \ $ Then (see e.g. \cite{CL}, \cite{S2})

\begin{equation}
\label{woodland}
\sum^\infty_{n=1}\frac{C_n(L)}{n!}z^n=\sum^\infty_{n=1}\frac{V_n(L)}{n!}
(e^z-1)^n.
\end{equation}

To apply the machinery of the cluster functions to the problem at hand, we consider a so-called $s$-modified random point process, which can
be constructed in the following way. We start with a random configuration $\xi=(x_i)_{i=-\infty}^{+\infty}$ from the
original random point field, and keep only those points $x_i$ for which there is exactly one neighbor to the right within distance $s$, 
i.e. $ x_{i+1}-x_i \leq s, \ \ x_{i+2}-x_{i}>s. \ $  The points $x_i$
for which this conditions is not satisfied are thrown away. As a result, we obtain a new random configuration $\xi(s) \subset \xi, $ such that
$\xi(s)= \{x_i : x_{i+1}-x_i \leq s, \ \ x_{i+2}-x_{i}>s \}, \ $
where $ \ldots < x_{-2} < x_{-1} < x_0 <x_1 <x_2<\ldots $ are the points of the original 
configuration $\xi. \ $  
The number of spacings  less than $s$ of the original random point field 
in the interval $ [0,L] $ is related 
to the number of points of the $s$-modified random point field in $[0,L].\ $
As will be shown later, for large $L$ and $s \sim L^{-1/3}, \ $ these two counting random variables coincide with probability very close to 1.

Since the moments and the cumulants of the counting random variable
$\#([0,L])$ can be expressed in terms of the integrals of point correlation and cluster functions (\ref{tel}), (\ref{woodland}), it is 
essential to be
able to calculate the point correlation and cluster functions of the $s$-modified random point process. We shall denote the 
$k$-
point correlation and $k$-point cluster functions of the modified random process by $\rho_k(x_1, \ldots, x_k; s)$ and 
$r_k(x_1, \ldots, x_k;s)$,
correspondingly. It follows from the inclusion-exclusion principle that,
provided $|x_i-x_j|>s, \ \ 1\leq i\neq j \leq k, \ $ one obtains

\begin{eqnarray}
\nonumber
& & \rho_k(x_1, \ldots, x_k;s) = 
\sum_{m=0}^{+\infty} \frac{(-1)^m}{m!} \* \int_{x_1}^{x_1+s} \ldots \int_{x_k}^{x_k+s} \* 
\int_{I(x_1,\ldots, x_k;s)^m}  \\
\label{sacto}
& & \rho_{2\*k+m}(x_1,\ldots, x_k, y_1, \ldots, y_k, z_1, \ldots, z_m) \*
dy_1 \ldots dy_k \* dz_1 \* \ldots dz_m,
\end{eqnarray}
where $I(x_1, \ldots, x_k;s)= \bigcup_{i=1}^m [x_i, x_i+s], $  and $ \ I(x_1,\ldots, x_k;s)^m=
I(x_1,\ldots, x_k;s) \times \ldots   \times I(x_1,\ldots, x_k;s) \  $
stands for the $m$-th fold Cartesian product of $I(x_1, \ldots, x_k;s)$
(see e.g. \cite{S1}, \cite{S2}).

In the determinantal case (\ref{detcorrfun}) the formula for the $k$-point cluster function of the $s$-modified random process has a 
somewhat similar structure (\cite{S1}, \cite{S2}). Provided
$|x_i-x_j|>s, \ \ 1\leq i\neq j \leq k, \ $ one has

\begin{eqnarray}
\nonumber
& & r_k(x_1, \ldots, x_k;s) = 
\sum_{m=0}^{+\infty} \frac{(-1)^m}{m!} \* \int_{x_1}^{x_1+s} \ldots \int_{x_k}^{x_k+s} \* 
\int_{I(x_1,\ldots, x_k;s)^m} \\
\label{westsacto}
& & \rho^{trun}_{2\*k+m}(x_1,\ldots, x_k, y_1, \ldots, y_k, z_1, \ldots, z_m) \*
dy_1 \ldots dy_k \* dz_1 \* \ldots dz_m,
\end{eqnarray}
where 
\begin{eqnarray}
\nonumber
& & \rho_{2\*k +m}^{trun}(x_1, \ldots, x_k, y_1, \ldots, y_k, z_1, \ldots, z_m)= 
\sum_{\sigma \in S_{2\*k +m}}^*
(-1)^{\sigma} \*  K(x_1,\sigma (x_1))\ldots K(x_k, \sigma(x_k)) \times \\
\label{eastsacto}
& &
 K(y_1, \sigma(y_1)) \ldots K(y_k, \sigma (y_k)) \* K(z_1, \sigma(z_1))
\ldots K(z_m,\sigma (z_m)),
\end{eqnarray}
where the summation in (\ref{eastsacto}) is over the permutations $\sigma \in S_{2\*k +m} $
satisfying
the property {\bf A}  described below (we note that $\sigma$ acts on the set of $ 2\*k +m \ $ variables $ \ x_1,\dots ,x_k, y_1,
\dots, y_k, z_1,\dots z_m $).

{\bf Property A}

{\it Consider an index $1 \leq i \leq k.\ $  Define $X(i)$ to be the subset of the set of variables $\{x_1,\dots ,x_k, y_1,
\dots, y_k, z_1,\dots z_m\}, \ $ that consists of $x_i, \  y_i, \ $ and those of the variables $z_1, \ldots, z_m, $ that belong to 
$[x_i, x_i +s] \ $ (we remind the reader that each $z_l, \ 1\leq l \leq m, \ $  
belongs to exactly one interval $ [x_j, x_j +s], \ 
1 \leq j \leq k). \ $  Then for any pair of disjoint integers $ 1 \leq i \neq j \leq k \ $ there exists a positive integer
$r=r(i,j) >0, \ $ such that $ \sigma^r (X(i))\bigcap X(j) \neq \emptyset.$}

We would like to bring to the reader's attention the fact that the relation between (\ref{sacto}) and (\ref{westsacto}) is, in a sense, 
quite similar to the relation between (\ref{winters}) and (\ref{davis}).

\section{Proof of the Main Result}
\label{sec:3}

The strategy of the proof is the following. We consider the rescaling  $\tilde{s}=s\*L^{-\frac{1}{3}}, \ $  
(we shall show that the smallest spacings in the interval 
$[0, L] $ 
are of order  $L^{-1/3}).$  We shall keep $s$ fixed as $ L \to \infty, \ $ so $\tilde{s} $ will be proportional to
$L^{-\frac{1}{3}}.\ $ We are interested in the asymptotics of
the integrals
\begin{equation}
\label{detroit}
V_k(L)= \int_{[-L,L]^k} r_k(x_1,x_2,\ldots x_k; \tilde{s}) \* dx_1 \* dx_2 \ldots dx_k.
\end{equation}
We claim that
$ \lim_{L\to \infty} V_1(L) = \alpha \* s^3, \ $ where $\alpha $ is defined in (\ref{nj}), and
$ \lim_{L\to \infty} V_k(L) = 0, \ $ for $k>1. \ $ 

\begin{lemma}
Let $V_k(L)$ be defined as in (\ref{detroit}), where 
$ r_k(x_1,x_2,\ldots x_k; \tilde{s}) $ is the $k$-point cluster function of the $\tilde{s}$-modified random point process introduced above 
and $\tilde{s}= s \* L^{-\frac{1}{3}}. $ Then
\begin{eqnarray}
\label{lemmochka}
\lim_{L \to \infty} V_k(L)=
\left\{ \begin{array}{cc} \alpha \* s^3  & { if \ } k=1, \\ 
                                                   0 & { if \ } k>1,
                                       \end{array} \right.
\end{eqnarray}
where $\alpha$ has been  defined  in (\ref{nj}).
\end{lemma}

The result of Lemma 1, combined with (\ref{woodland}), implies that
the number of the points of the $\tilde{s}$- modified random process in the interval $[0,L]$ converges in distribution to the Poisson law
as $ L \to \infty.$ 

Once Lemma 1 is proven , we shall show that the number of points in $[0, L]$ of the original determinantal process that have at least 
two neighbors to the 
right within distance $s/L^{1/3}$ is zero with probability very close to 1, provided that $L$ is large and $s$ stays finite.

{\bf Proof of Lemma}

We start with $V_1(L).\ $ Consider the one-point correlation function (intensity) of the $\tilde{s}$-modified point process 
$\rho_1(x;\tilde{s}).\ $ 

\begin{equation}
\label{sanfran}
\rho_1(x;\tilde{s}) = 
\sum_{m=0}^{+\infty} \frac{(-1)^m}{m!} \* \int_x^{x+\tilde{s}} \ldots
\int_x^{x+\tilde{s}} \rho_{m+2}(x, y, z_1, \ldots, z_m) \* dy \* d^m z.
\end{equation}

We claim that in the determinantal case
\begin{equation}
\label{LA}
\rho_1(x;\tilde{s})= \int_x^{x+\tilde{s}} \rho_2(x,y)\* D(x,y;\tilde{s}) \* dy,
\end{equation}
where $ \ D(x,y;\tilde{s}) \ $ is the Fredholm determinant of the integral operator $\tilde{K}$ on $\ L^2([x, x+\tilde{s}]), \ $ 
\begin{equation}
\label{NYC}
D(x,y;\tilde{s})= \det(1-\tilde{K}), \ \ \ \tilde{K}:L^2([x,x+\tilde{s}])\to L^2([x,x+\tilde{s}]).
\end{equation} 
The kernel
$ \ \tilde{K}(u,v) \ $ in (\ref{NYC}), (\ref{LA}) depends on $x$ and $\  y, \ $  and is given by the formula
\begin{eqnarray}
\nonumber
\tilde{K}(u,v) &=& K(u,v)- K(u,x)\*T_{11}(x,y) \* K(x,v) -K(u,x)\*T_{12}(x,y) \* K(y,v) \\
\label{detroit}
  &-& K(u,y)\*T_{21}(x,y) \* K(x,v) -
K(u,y)\*T_{22}(x,y) \* K(y,v),
\end{eqnarray}
where 
\[ \left( \begin{array}{cc} T_{11}(x,y) &  T_{12}(x,y) \\
T_{21}(x,y) & T_{22}(x,y) \end{array} \right)=
\left( \begin{array}{cc} K(x,x) &  K(x,y) \\
K(y,x) & K(y,y) \end{array} \right)^{-1}. 
\]
Indeed,  let us introduce the notation $\ K[x_1, \ldots, x_k]:= \det(K(x_i, x_j))_{i,j=1, \ldots, k}. \ $ Then,
$$ K[x,y,z_1, \ldots, z_m]=K[x,y]\* \tilde{K}[z_1, \ldots, z_m]=\rho_2(x,y)\*\tilde{K}[z_1, \ldots, z_m].$$  
In other words, the conditional distribution of a determinantal random 
point process with the correlation kernel $K$, given there are two particles at $ x $ and $ y $ is again a determinantal random point 
process (on $\mathcal{R}^1 \setminus\{x,y \}) \ $ with the kernel $\ \tilde{K} \ $ (see e.g. \cite{ST2}). This allows us to rewrite 
(\ref{sanfran}) as
\begin{equation}
\label{sanfran5}
\rho_1(x;\tilde{s}) = 
\int_x^{x+\tilde{s}} \rho_2(x,y) \* \bigl(\sum_{m=0}^{+\infty} \frac{(-1)^m}{m!} \*
\int_{[x, x+\tilde{s}]^m} \tilde{K}(z_1, \ldots, z_m) \* d^m z \bigr) \*dy,
\end{equation}
and (\ref{LA}) follows.

The intensity
$\rho_1(x;\tilde{s}) \ $ is constant (i.e. it does not depend on $x$) in the translation-invariant case. 
To estimate $\rho_1(x;\tilde{s})=\rho_1(0;\tilde{s}), \ $ 
we note that
\begin{eqnarray}
\nonumber
& &\rho_{m+2}(x,y,z_1, \ldots, z_m)=K[x,y, z_1, \ldots, z_m] \leq K(x,x)\*K(y,y) \*K(z_1, z_1) \ldots K(z_m, z_m) \\
\label{hruhru}
& &\leq g(o)^{m+2}, 
\end{eqnarray}
since the determinant of a non-negative definite matrix is bounded from above by the product of the diagonal entries (for the generalization 
of this result see Lemma 2 below).
Then 

\begin{eqnarray}
\nonumber
& & \rho_1(x;\tilde{s}) = 
\int_x^{x+\tilde{s}} \rho_2(x,y) \* dy  - \int_x^{x+\tilde{s}} \int_x^{x+\tilde{s}} \rho_3(x,y,z_1) \* dy\* dz_1 + \\
\nonumber
& & \frac{1}{2} \*
\int_x^{x+\tilde{s}} \int_x^{x+\tilde{s}} \int_x^{x+\tilde{s}}
\rho_4(x,y,z_1, z_2) \* dy \* dz_1 \* dz_2 + \\
\nonumber
& & \sum_{m=3}^{+\infty} \frac{(-1)^m}{m!} \* \int_x^{x+\tilde{s}} \ldots
\int_x^{x+\tilde{s}} \rho_{m+2}(x, y, z_1, \ldots, z_m) \* dy \* d^m z = \\
\nonumber
& & \int_x^{x+\tilde{s}} \rho_2(x,y) \* dy  - \int_x^{x+\tilde{s}} \int_x^{x+\tilde{s}} \rho_3(x,y,z_1) \* dy\* dz_1 + \\
\nonumber
& & \frac{1}{2} \*
\int_x^{x+\tilde{s}} \int_x^{x+\tilde{s}} \int_x^{x+\tilde{s}}
\rho_4(x,y,z_1, z_2) \* dy \* dz_1 \* dz_2 + O(\tilde{s}^4).
\end{eqnarray}

To estimate $\ \rho_3(x,y,z_1)= K[x,y,z_1] \ $ and $ \ \rho_4(x,y,z_1, z_2)=K[x,y,z_1, z_2] \ $ we recall that the point 
correlation functions 
are given by the determinants,
and subtract the first column from the other columns both in $ \ K[x,y,z_1] \ $ and $ \ K[x,y,z_1, z_2]. \ $
Since $ y\in [x, x+\tilde{s}], \ z_i \in [x, x+\tilde{s}], \ i\geq 1, \ $ and the first derivative of $g$ is uniformly bounded, 
we observe that 
$\rho_3(x,y,z_1) =K[x,y,z_1]=O(\tilde{s}^2), \ \ \rho_4(x,y,z_1, z_2)=
K[x,y,z_1, z_2]=O(\tilde{s}^3), \ $ and, therefore
\begin{equation}
\label{princeton}
\rho_1(x;\tilde{s})= \int_x^{x+\tilde{s}} \rho_2(x,y) \* dy \ + O(\tilde{s}^4) = \int_0^{\tilde{s}} (g^2(0) - g^2(t)) \* dt + 
O(\tilde{s}^4) = \alpha \* \tilde{s}^3 + O(\tilde{s}^4), 
\end{equation}
where $ \ \alpha \ $ has been defined in (\ref{nj}).
It follows from (\ref{princeton}) that 
$ \ \lim_{L \to \infty} V_1(L) = \lim_{L \to \infty} \int_0^L \rho_1(x;\tilde{s}) \* dx=\lim_{L \to \infty} \rho_1(0;\tilde{s})\*L=
\alpha \*s^3. \ $ 

Next, we show that $ \ \lim_{L \to \infty} V_k(L) =0 \ $ for $ \ k >1. \ $ We remind the reader that $V_k(L)$ has been defined as
$ \ V_k(L)= 
\int_{[-L,L]^k} r_k(x_1,x_2,\ldots x_k; \tilde{s}) \* dx_1 \* dx_2 \ldots dx_k. \ $
We start with the case $ k=2. $  Recall (see (\ref{westsacto})) that for $ \ |x_1-x_2| >\tilde{s} \ $
\begin{eqnarray}
\nonumber
& & r_2(x_1, x_2;\tilde{s}) = 
\sum_{m=0}^{+\infty} \frac{(-1)^m}{m!} \* \int_{x_1}^{x_1+\tilde{s}} \int_{x_2}^{x_2+\tilde{s}} \* 
\int_{I(x_1, x_2;\tilde{s})^m}
\rho^{trun}_{4+m}(x_1, x_2, y_1, y_2, z_1, \ldots, z_m) \\
\label{westsacto2}
& & dy_1 \* dy_2 \* dz_1 \* \ldots dz_m,
\end{eqnarray}
where $ \ \rho_{4+m}^{trun}(x_1, x_2, y_1, y_2, z_1, \ldots, z_m) \ $ has been defined in (\ref{eastsacto}).

As described in the Property A (right after the formula (\ref{eastsacto})), in order to define 
$ \ \rho_{4+m}^{trun}(x_1, x_2, y_1, y_2, z_1, \ldots, z_m) \ $ one introduces a partition
$ \ X(1) \bigsqcup X(2) = 
\{x_1, x_2, y_1, y_2, z_1, \ldots, z_m\}, \ $
where $ \ X(1) \ $ consists of  $ \ x_1, y_1, $ and those of the variables $ \ z_1, \ldots, z_m \ $ that belong to $ \ [x_1, x_1 +s], \ $
and $ \ X(2) \ $ consists of  $ \ x_2, y_2, $ and those of the variables $ \ z_1, \ldots, z_m \ $ that belong to $ \ [x_2, x_2 +s]. \ $
Let $ \ X(1) \bigcap \{ z_1, \ldots, z_m \} = \{z_{i_1}, \ldots, z_{i_l} \}, \ $ and $ \ X(2) \bigcap \{ z_1, \ldots, z_m \} 
= \{z_{j_1}, \ldots, z_{j_{m-l}} \}. \ $ Then 
\begin{eqnarray}
\nonumber
& & \rho_{4+m}^{trun}(x_1, x_2, y_1, y_2, z_1, \ldots, z_m) = K[x_1, y_1, x_2, y_2, z_1, \ldots, z_m] - \\
\label{folsom}
& & K[x_1, y_1, z_{i_1}, \ldots, z_{i_l}]\*
K[x_2, y_2, z_{j_1}, \ldots, z_{j_{m-l}}]. 
\end{eqnarray}

We claim that 
\begin{equation}
\label{yuba}
|\rho_{4+m}^{trun}(x_1, x_2, y_1, y_2, z_1, \ldots, z_m)| \leq (m+4)! \* (C\*\tilde{s})^{2+m} \* 
\frac{const}{1+|x_1-x_2|^{1+ \epsilon}}, 
\end{equation}
where $ const $ is a constant that may depend on $s,$ 
and $C$ is the constant introduced after the formulas (\ref{dva}), (\ref{tri}).

The factor $ \ (C \*\tilde{s})^{2+m} \ $ in (\ref{folsom}) follows from the uniform bound on the derivative of $g,$ and the fact that the
$m+2$ variables $y_1, y_2, z_1, \ldots, z_m$ are within distance $\tilde{s}$ from either $x_1$ or $x_2. \ $ In other words, one can subtract 
the first 
column in the matrices in $ K[x_1, y_1, z_{i_1}, \ldots, z_{i_l}] $ and $ K[x_2, y_2, z_{j_1}, \ldots, z_{j_{m-l}}] $ from the 
other columns,  and subtract the first and the third column in the matrix in $ K[x_1, y_1, x_2, y_2, z_1, \ldots, z_m] $ from the 
corresponding columns. Such linear operations do not change the value of the determinants, and the new matrices will contain the terms
$ g(u-w)-g(x_j-w),  $ in all columns, except those corresponding to $x_1$ and $x_2, \ $  where $ j=1,2, $ and $ u \in [x_j, x_j+s]. $ Such 
terms can be estimated from above by $ \bigl(\max_{x \in[x_j-w, x_j+\tilde{s}-w]} |g'(x)|\bigr) \* \tilde{s}. \ $  

It follows from the definition that $ \rho_{4+m}^{trun}(x_1, x_2, y_1, y_2, z_1, \ldots, z_m) \ $ can be 
written as a sum over at most $ (m+4)! $ permutations, each term being a product $m+4$ factors. As we just showed, 
$m+2$ out of those $m+4$ factors can be 
estimated in absolute value by $ C \* \tilde{s}. \ $ 
Moreover, Property A implies that  at least two factors in each term  must be given either by
$ g(x_1-v), \ $ or by $ \ g(u -v) - g(x_1 - v), \ $ or by $g(x_2-u)$, or by
$ g(v-u)- g(x_2-u), \ $ where $ \ u \in [x_1, x_1 +\tilde{s}], \ v \in [x_2, x_2 +\tilde{s}]. \ $ 
The inequalities (\ref{dva}), (\ref{tri})  imply that
these two factors each contribute an upper bound $ \frac{C}{1 + (|x_1-x_2|-\tilde{s})^{\frac{1}{2}+\epsilon}} $ and the desired estimate
(\ref{yuba}) follows. 

We recall that we defined above $ I(x_1, x_2;\tilde{s})= [x_1, x_1 +\tilde{s}] 
\bigcup [x_2, x_2+\tilde{s}]. $ Then

\begin{eqnarray}
\nonumber
& & |\int_{x_1}^{x_1+\tilde{s}} \int_{x_2}^{x_2+\tilde{s}} \* 
\int_{I(x_1, x_2;\tilde{s})^m}
\rho^{trun}_{4+m}(x_1, x_2, y_1, y_2, z_1, \ldots, z_m)
dy_1 \* dy_2 \* dz_1 \* \ldots dz_m | \leq \\
\label{kengurchik}
& & (m+4)! \* (C\*\tilde{s})^{2+m} \* (2\*\tilde{s})^{2+m} \* \frac{const}{1+|x_1-x_2|^{1+\epsilon}},
\end{eqnarray}
and
\begin{eqnarray}
\nonumber
& & |r_2(x_1, x_2;\tilde{s})| \leq  
\sum_{m=0}^{+\infty} \frac{1}{m!} \* (m+4)! \* (C\*\tilde{s})^{2+m} \* (2\*\tilde{s})^{2+m} \*
\frac{const}{1+|x_1-x_2|^{1+\epsilon}} \leq \\
\label{sandiego}
& & \frac{const\* (\tilde{s})^4}{1+|x_1-x_2|^{1+\epsilon}}.
\end{eqnarray}
We remind the reader that (\ref{sandiego}) has been derived for $ \ |x_1-x_2|>\tilde{s}. $ 

Since $\tilde{s}=s\*L^{-1/3}, \ $  it follows from the above estimate that
\begin{equation}
\label{fairfield}
\lim_{L \to \infty} \* \int_0^L \int_0^L \ r_2(x_1, x_2; \tilde{s}) \* \chi_D(x_1, x_2) \* dx_1 \* dx_2 =\lim_{L \to \infty} 
O(\tilde{s}^4\*L) =0,
\end{equation}
where $ D=\{ (x_1, x_2): \ |x_1-x_2| > s\*L^{-1/3} \}.$

To estimate $r_2(x_1, x_2; \tilde{s}) $ on $D^c,$ we note that $$ r_2(x_1, x_2; \tilde{s})=\rho_2(x_1, x_2; \tilde{s}) - 
\rho_1(x_1;\tilde{s}) \* \rho_1(x_2; \tilde{s}).$$

It follows then from (\ref{princeton}) that $ \rho_1(x_1;\tilde{s})= \rho_1(x_2;\tilde{s}) \leq const \* s^3 \* L^{-1}. \ $
To estimate $\rho_2(x_1, x_2; \tilde{s})$ we can assume without  loss of generality that  $ x_2 \leq x_1 \leq x_2+\tilde{s} $. Then
\begin{eqnarray}
\nonumber
& & \rho_2(x_1, x_2; \tilde{s}) = 
\sum_{m=0}^{+\infty} \frac{(-1)^m}{m!} \* \int_{x_1}^{x_1+s} \* 
\int_{I(x_1,x_2;\tilde{s})^m}
\rho_{3+m}(x_1,x_2, y_1, z_1, \ldots, z_m) \\
\label{westsacto11}
& & dy_1 \* dz_1 \* \ldots dz_m,
\end{eqnarray}
where $ I(x_1, x_2;\tilde{s})= [x_1, x_1+\tilde{s}] \bigcup [x_2, x_2+\tilde{s}]= [x_2, x_1+\tilde{s}].$
Subtracting the first column in $ \ \rho_{3+m}(x_1,x_2, y_1, z_1, \ldots, z_m)= K[x_1,x_2, y_1, z_1, \ldots, z_m] \ $ from the other columns
and using (\ref{tri}), we see that $ \ \rho_{3+m}(x_1,x_2, y_1, z_1, \ldots, z_m) \leq (m+3)! \* (const*\tilde{s})^{m+2}. $ Integrating 
over $  y_1, z_1, 
\ldots, z_m $ and summing over $m$ we obtain 
\begin{equation}
\label{lodi}
 \rho_2(x_1, x_2; \tilde{s}) \leq const \* \tilde{s}^3,
\end{equation}
which implies

\begin{eqnarray}
\nonumber
& & \lim_{L \to \infty} \* \int_0^L \int_0^L \ r_2(x_1, x_2; \tilde{s}) \* \chi_{D^c}(x_1, x_2) \* dx_1 \* dx_2= \\
\nonumber
& & \lim_{L \to \infty} \* \int_0^L \int_0^L \ \rho_2(x_1, x_2; \tilde{s}) \* \chi_{D^c}(x_1, x_2) \* dx_1 \* dx_2 - \\
\nonumber
& & \lim_{L \to \infty} \int_0^L \int_0^L \ \rho_1(x_1; \tilde{s})\*\rho_1(x_1; \tilde{s}) \*\chi_{D^c}(x_1, x_2)\* dx_1 \* dx_2=\\
\label{martinez}
& &
\lim_{L \to \infty} \bigl(O(\tilde{s}^4)\*L
- O(\tilde{s}^7)\*L\bigr)=0.
\end{eqnarray}

Combining (\ref{fairfield}) and (\ref{martinez}) one obtaines
$
\lim_{L \to \infty} V_2(L)=0.
$

The argument in the case of general $k>2$ is quite similar. Again,
we first estimate $ r_k(x_1, x_2, \ldots, x_k;\tilde{s}) $ on $D=\{ (x_1, x_2, \ldots, x_k): |x_i-x_j|>\tilde{s}, \ 1\leq i\neq j 
\leq k\}. \  $
We will use formulas (\ref{westsacto}) and (\ref{eastsacto}).  To estimate 
$ \rho^{trun}_{2k+m}(x_1, \ldots, x_k, y_1, \ldots, y_k, z_1, \ldots, z_m) $ we consider the
partition $ X(1) \bigsqcup \ldots \bigsqcup X(k)= \{x_1, \ldots, x_k, y_1, \ldots, y_k, z_1, \ldots, z_m \}, $
where 
$\ X(i)= $ \\ $ \{x_1, \ldots, x_k, y_1, \ldots, y_k, z_1, \ldots, z_m \} \bigcap [x_i, x_i+\tilde{s}]. \ $ Let $ X(i) \bigcap
\{ z_1, \ldots, z_m \} = \{z_1^{(i)}, \ldots, z_{n_i}^{(i)} \}. \ $ Then
$X(i)= \{x_i, y_i, z_1^{(i)}, \ldots, z_{n_i}^{(i)} \}. $

It follows from Property A and the inclusion-exclusion principle that
\begin{equation}
\label{baker}
\rho^{trun}_{2k+m}(x_1, \ldots, x_k, y_1, \ldots, y_k, z_1, \ldots, z_m) = \sum_{G} (-1)^j \* (j-1)! \prod_{l=1}^j K_l,
\end{equation}
where the summation is over all partitions $G=G_1 \bigsqcup \ldots \bigsqcup G_j $
of $[k]=\{1,2,\ldots, k\},\  \ j=1, \ldots, k, $ and $
K_l= K[x_i, y_i, z_1^{(i)}, \ldots, z_{n_i}^{(i)} : i \in G_l ]; \ $ in other words, $K_l$ depends on the variables from
$ \ \bigsqcup_{i \in G_l} X(i), \ $ and it is given by the determinant of the matrix built from the correlation kernel $K(x,y).\ $
We claim that 
\begin{eqnarray}
\nonumber
& &|\rho^{trun}_{2k+m}(x_1, \ldots, x_k, y_1, \ldots, y_k, z_1, \ldots, z_m)| \leq (m+2k)! \* (C\*\tilde{s})^{k+m} \times \\
\label{esparto1}
& &\bigl(
\frac{const}{1+|x_1-x_2|^{\frac{1}{2}+\epsilon}} \* \frac{const}{1+|x_2-x_3|^{\frac{1}{2}+\epsilon}}\cdots
\frac{const}{1+|x_k-x_1|^{\frac{1}{2}+\epsilon}}  + \ ... \ 
\bigr),
\end{eqnarray}
where the summation in the last factor of the r.h.s. of (\ref{esparto1}) is over all $(k-1)!$ cyclic permutations 
(for example, the first term in the sum corresponds to the cyclic permutation
$ 1 \to 2\to 3 \to \ldots \to k \to 1).$ 
We claim that the estimate (\ref{esparto1}) follows from (\ref{dva}), (\ref{tri}), (\ref{baker}) and Property A. 
As in the case $k=2$ discussed above, we use the fact that each of the $k+m$ variables $y_1, \ldots, y_k, z_1, \ldots, z_m $ lies 
within distance
$ \tilde{s} $ from one of the $x_i$'s, $ i=1, \ldots, k.$ In each $K_l=  K[x_i, y_i, z_1^{(i)}, \ldots, z_{n_i}^{(i)} : i \in G_l ] $ in
(\ref{baker}) we subtract for each $ i \in G_l \ $ the column corresponding to $x_i $ from the column corresponding to $y_i$ and from the 
other columns corresponding 
to the variables from $X(i).$ These linear operations do not change the values of determinants, and, therefore, do not change the value of
$\rho^{trun}_{2k+m}(x_1, \ldots, x_k, y_1, \ldots, y_k, z_1, \ldots, z_m).$ Now, according to the Property A,  we observe that 
$\rho^{trun}_{2k+m} $ is a sum 
of at most $(m+2k)! $ terms.   Each term is a product of $m+2k$ factors. 
Property A  assures that  each term in the sum  can be put into correspondence with
a cyclic permutation $\sigma $ on the set of $k$ variables $x_1, x_2, \ldots, x_k, $ in such a way that  $k$ out of $m+2k$ terms in the 
product are of the form
$ g(x_{\sigma(i)}-v), \ $ or $ \ g(u -v) - g(x_{\sigma(i)} - v), \ i=1, \ldots, k, \ \ $ where $
v \in [x_{\sigma(i+1)}, x_{\sigma(i+1)} +\tilde{s}], \ $ and $ \sigma (k+1)=\sigma (1). \ $  The bounds (\ref{dva}), (\ref{tri}) then imply
(\ref{esparto1}) in the same manner as has been shown in the case $k=2.\ $
Therefore,
\begin{eqnarray}
\nonumber
 &|&\int_{x_1}^{x_1+\tilde{s}} \ldots \int_{x_k}^{x_k+\tilde{s}} \* \int_{I(x_1, \ldots, x_k; \tilde{s})^m}
\rho^{trun}_{2k+m}(x_1, \ldots, x_k, y_1, \ldots, y_k, z_1, \ldots, z_m) \* dy_1 \ldots dy_k \\
\nonumber
&\*& dz_1 \ldots dz_m| 
\  \leq   (m+2k)! \* (C\*\tilde{s})^{k+m} \* (k\*s)^{k+m} \times \\
\label{rocklin}
& & \bigl(
\frac{const}{1+|x_1-x_2|^{\frac{1}{2}+\epsilon}} \* \frac{const}{1+|x_2-x_3|^{\frac{1}{2}+\epsilon}}\cdots
\frac{const}{1+|x_k-x_1|^{\frac{1}{2}+\epsilon}}  + ...
\bigr),
\end{eqnarray}
provided $ \vec{x}=(x_1, \ldots, x_k) \in D, \ $ i.e. $|x_i-x_j|>\tilde{s} $ for $i\neq j.$
Then on $D$ we have an estimate
\begin{eqnarray}
\nonumber
& & |r_k(x_1, \ldots, x_k;\tilde{s})| \leq  
\sum_{m=0}^{+\infty} \frac{1}{m!} \* (2\*k+m)! \* (C\*\tilde{s})^{k+m} \* (k\*\tilde{s})^{k+m} \times \\
\nonumber
& & \bigl(
\frac{const}{1+|x_1-x_2|^{\frac{1}{2}+\epsilon}} \times \frac{const}{1+|x_2-x_3|^{\frac{1}{2}+\epsilon}}\times \cdots
\frac{const}{1+|x_k-x_1|^{\frac{1}{2}+\epsilon}}  + ...
\bigr) \leq \\
\nonumber
& & Const_k \* (\tilde{s})^{2\*k} \* \bigl(
\frac{const}{1+|x_1-x_2|^{\frac{1}{2}+\epsilon}} \times \frac{const}{1+|x_2-x_3|^{\frac{1}{2}+\epsilon}}\times\cdots
\frac{const}{1+|x_k-x_1|^{\frac{1}{2}+\epsilon}}  + ...
\bigr),
\end{eqnarray}
and
\begin{equation}
\nonumber
|\int_0^L \ldots \int_0^L \ r_k(x_1,\ldots, x_k; \tilde{s}) \* \chi_D(x_1,\ldots, x_k) \* dx_1 \ldots dx_k| \leq 
Const_k \* (L^{1+(k-1)\*(\frac{1}{2}-\epsilon)}) \* s^{2 \*k} \* L^{-\frac{2\*k}{3}}.
\end{equation}
The last estimate implies
\begin{equation}
\label{imperial}
\lim_{L \to \infty} \int_0^L \ldots \int_0^L \ r_k(x_1,\ldots, x_k; \tilde{s}) \* \chi_D(x_1,\ldots, x_k) \* dx_1 \ldots dx_k =0,
\end{equation}
for $k \geq 3. \ $
Our next goal is to show that
\begin{equation}
\label{hockey}
\lim_{L \to \infty} \int_0^L \ldots \int_0^L \ r_k(x_1,\ldots, x_k; \tilde{s}) \* \chi_D^c(x_1,\ldots, x_k) \* dx_1 \ldots dx_k =0.
\end{equation}
To estimate $ r_k(x_1,\ldots, x_k; \tilde{s})$ on $D^c,$ we rewrite the formula (\ref{cluster}) that expresses the $k$-point cluster function 
in terms of point correlation functions:
\begin{eqnarray}
\nonumber
& & r_k(x_1,\ldots, x_k; \tilde{s})= \rho_k(x_1,\ldots, x_k; \tilde{s})-\rho_1(x_1; \tilde{s}) \* 
\rho_{k-1}(x_2, x_3, \ldots, x_k; \tilde{s}) - \\
\nonumber
& &\rho_1(x_2; \tilde{s}) \* \rho_{k-1}(x_1, x_3, \ldots, x_k; \tilde{s}) \ldots
-\rho_1(x_k; \tilde{s}) \* \rho_{k-1}(x_1, x_2, \ldots, x_{k-1}; \tilde{s}) +\\
\nonumber
& & 2\*\rho_2(x_1, x_2; \tilde{s}) \* \rho_{k-2}(x_3, \ldots, x_k; \tilde{s}) +
2 \*\rho_2(x_1, x_3; \tilde{s}) \* \rho_{k-2}(x_2, x_4, \ldots, x_k; \tilde{s}) + \\
\label{oxford}
& & \ldots \rho_2(x_{k-1}, x_k; \tilde{s}) \* \rho_{k-2}(x_1, x_2, \ldots, x_{k-2}; \tilde{s}) - \ldots
\end{eqnarray}
We claim that the integral of each of the terms in (\ref{oxford}) over $ [0,L]^m \bigcap D^c $ has a zero limit as $ L \to \infty. $
To prove it, we consider an arbitrary term in (\ref{oxford}),
\begin{equation}
\label{dubna}
\rho_{k_1}(x_1,x_2, \ldots, x_{k_1}; \tilde{s}) \* \rho_{k_2}(x_{k_1+1}, x_{k_1+2}, \ldots, x_{k_1+k_2}; \tilde{s})\cdots
\rho_{k_l}(x_{k_1+\ldots k_{l-1}+1},\ldots, x_k; \tilde{s}), 
\end{equation}
where $ k_1+k_2 \ldots + k_l=k, \ k_i \geq 1, \ i=1, \ldots, k. \ $
We shall estimate $ \rho_{k_1}(x_1,x_2, \ldots, x_{k_1}; \tilde{s}), $ the other $ l-1 $ factors are estimated in the same way.

First assume that none of the variables $ x_1, x_2, \ldots, x_{k_1} $ are within distance $\tilde{s} $ from each other.
Then one can clearly estimate $ \rho_{k_1}(x_1,x_2, \ldots, x_{k_1}; \tilde{s}) $ from above as
\begin{equation}
\label{khanin}
\rho_{k_1}(x_1,x_2, \ldots, x_{k_1}; \tilde{s}) \leq \int_{x_1}^{x_1+\tilde{s}} \ldots \int_{x_{k_1}}^{x_{k_1}+\tilde{s}} \rho_{2\*k_1} 
(x_1, \ldots, x_{k_1},
y_1, \ldots y_{k_1}) \* dy_1\ldots dy_{k_1}.
\end{equation}
Now, since $ \rho_{2\*k_1} (x_1, \ldots, x_{k_1},
y_1, \ldots y_{k_1}) = K[x_1, \ldots, x_{k_1}, y_1, \ldots, y_{k_1}], $ and \\ 
$ K[x_1, \ldots, x_{k_1}, y_1, \ldots, y_{k_1}] $ 
is the determinant of a 
$ (2\*k_1)$-dimensional (non-negative definite) real symmetric matrix, we can estimate the determinant from above by the product of the 
determinants
of the $ \ 2\times 2 $ diagonal blocks
\begin{equation}
\label{poezd}
K[x_1, \ldots, x_{k_1}, y_1, \ldots, y_{k_1}]=K[x_1,y_1, \ldots, x_{k_1}, y_{k_1}] \leq \prod_{i=1}^{k_1} K[x_i, y_i]. 
\end{equation}
The bound (\ref{poezd}) follows from the Fischer inequality which we state below as Lemma 2.

\begin{lemma}
Let $ M= \left(
\begin{array}{cc} A   & B \\
                                 B^* & C 
               \end{array} \right) $ be a 
block matrix, let $A $ and $C$ be $ n\times n $ and, respectively, $m \times m $ non-negative definite matrices, and $B$ be a
$m \times n $ matrix. Then

\begin{eqnarray}
\label{kuryane}
\det\left(
\begin{array}{cc} A   & B \\
                                 B^* & C 
               \end{array} \right) \leq \det A \* \det C.
\end{eqnarray}
\end{lemma}

{\bf Proof}

To prove Lemma 2, it is enough to reduce it to the obvious case
$M=\left(
\begin{array}{cc} Id & B \\
                                 B^* & Id
               \end{array} \right) $ by appropriate rotations and dilations in $C^n $ and $C^m$ (see e.g. \cite{ST2}).

As was shown above (see calculations around formula (\ref{princeton})
\begin{equation}
\label{glendale}
\int_{x_i}^{x_i+\tilde{s}} K[x_i, y_i]\* dy_i \leq const\* (\tilde{s})^3,
\end{equation}
which then implies that 
\begin{equation}
\label{pasadena1}
\rho_{k_1}(x_1,x_2, \ldots, x_{k_1}) \leq const^{k_1}\* (\tilde{s})^{3\* k_1}. 
\end{equation}

If none of the variables are within $\tilde{s}$ from each other in all factors  in (\ref{dubna}), We infer from (\ref{pasadena1}) that
\begin{eqnarray}
\nonumber
& & \rho_{k_1}(x_1,x_2, \ldots, x_{k_1}; \tilde{s}) \* \rho_{k_2}(x_{k_1+1}, x_{k_1+2}, \ldots, x_{k_1+k_2}; \tilde{s})\times \ldots \\
\label{straus}
& & \rho_{k_l}(x_{k_1+\ldots k_{l-1}+1},\ldots, x_k; \tilde{s}) 
\leq Const \* (\tilde{s})^{3\*k}= O(L^{-k}),
\end{eqnarray}
and the integral of
the l.h.s. of (\ref{straus}) over $ [0,L]^k \bigcap D^c $  goes to zero as $ L \to \infty,  $ since $ vol([0,L]^k \bigcap D^c)=O(L^{k-1}). $

If some of the variables in $\rho_{k_1}(x_1, \ldots, x_{k_1}) $  are within the distance $\tilde{s}$ from one another, the analysis is quite 
similar.  Let us assume, for example, that $ x_1 \leq x_2 \leq \ldots \leq x_{k_1}, $ and that 
$ x_i \leq x_{i+1}\leq x_{i}+\tilde{s}, \ \ i=1, 
\ldots, p,  \ $ and that the rest of the variables $ x_{p+1}, \ldots, x_{k_1} $ are not within the distance $\tilde{s} $ from each other.
Then
\begin{equation}
\label{khanin1}
\rho_{k_1}(x_1,x_2, \ldots, x_{k_1}; \tilde{s}) \leq \int_{x_{p+1}}^{x_{p+1}+\tilde{s}} \ldots \int_{x_{k_1}}^{x_{k_1}+\tilde{s}} 
\rho_{2\*k_1-p} (x_1, \ldots, x_{k_1},
y_{p+1}, \ldots y_{k_1}) \* dy_{p+1}\ldots dy_{k_1}.
\end{equation}
One can write
\begin{eqnarray}
\nonumber
& &\rho_{2\*k_1-p} (x_1, \ldots, x_{k_1},
y_{p+1}, \ldots y_{k_1}) = K[x_1, \ldots, x_{k_1}, y_{p+1}, \ldots, y_{k_1}] \leq K[x_1, x_2, \ldots, x_{p+1}, y_{p+1}]\\
\label{sosedi}
& &
\times \prod_{i=p+2}^{k_1}
K[x_i, y_i].
\end{eqnarray}
As before, 
\begin{equation}
\label{glendale1}
\int_{x_i}^{x_i+\tilde{s}} K[x_i, y_i]\* dy_i \leq const\* (\tilde{s})^3, \ \ i=p+1, \ldots, k_1.
\end{equation}
As for the term  $ K[x_1, \ldots, x_{p+1}, y_{p+1}], $ one can substract the first column from all other columns, and obtain
\begin{equation}
\label{uraura}
K[x_1, \ldots, x_{p+1}, y_{p+1}] \leq Const_k \* (\tilde{s})^{2p+2}, 
\end{equation}
since $|g(x)-g(y)|= O(\tilde{s^2}) \ $ for 
$ 0 \leq x, y \leq \tilde{s} \ $ (we used the fact that $ g'(0)=0 $).
Combining (\ref{glendale1}) and (\ref{uraura}), and integrating over the $y$'s we obtain
\begin{equation}
\label{elk1}
\rho_{k_1}(x_1,x_2, \ldots, x_{k_1}; \tilde{s})\leq Const \* (\tilde{s})^{3\*k_1-p+2}= O(L^{-k_1 +\frac{p}{3}- \frac{2}{3}}).
\end{equation}
Note, however, that
\begin{equation}
\label{elk2}
Vol \{ (x_1, \ldots, x_{k_1}) :  x_i \leq x_{i+1}\leq x_{i}+\tilde{s}, \ \ i=1, 
\ldots, p \} \bigcap [0,L]^{k_1}= O(L^{k_1-p}\*\tilde{s}^p),
\end{equation}
and the product of the right hand sides of (\ref{elk1}) and (\ref{elk2}) goes to zero.

If there are several factors in (\ref{straus}) for which there are variables within distance $\tilde{s} $ 
from each other, the analysis is very similar, and we leave the details to the reader.
Combining all the estimate together, one concludes the integral of  the l.h.s. of (\ref{straus}) over $ [0,L]^m \bigcap D^c $  goes 
to zero as $ \L \to \infty. $  This finishes the proof of Lemma.

The result of Lemma 1 and formula (\ref{woodland}) imply that  
$ \lim_{L \to\infty} \* \sum_{n=1}^{+\infty} \frac{C_n(L)}{n!} \* z^n = \alpha \* s^3 \* (e^z-1), $
where $ \{ C_n(L) \}_{n=1}^{+\infty} $ is the sequence of the cumulants of the counting random variable $N(L), $  where
$N(L)$  is the number of the points of the $\tilde{s}$-modified 
random pont field in the interval $[0, L].$ It follows from the definition of the $\tilde{s}$-modified 
random pont field that $N(L)=N_1(L)+N_2(L), \ $ where $ N_1(L) $ counts the number of particles of the original random point field that 
have exactly one neighbor within distance $\tilde{s} $ to the right, and $ N_2(L) $ counts the number of particles of the original random 
point field that 
have more than one neighbor within distance $\tilde{s} $ to the right. We claim that the probability that $N_2(L)\neq 0 $ is going to zero as 
$ L \to \infty. $ Specifically, we establish

\begin{lemma}
\begin{equation}
\label{skoro}
\lim_{L \to \infty} E \* N_2(L) =0,
\end{equation}
where $E$ denotes the mathematical expectation.
\end{lemma}

Since $N_2(L) $ is a non-negative, integer-valued random variable, (\ref{skoro}) implies that $ \Pr (N_2(L)\neq 0) \to 0  \ $ as 
$ L \to \infty. $

The proof of  Lemma 3 is elementary. We use the estimate
\begin{equation}
\label{kir}
E \* N_2(L) \leq \int_0^L \left( \int_x^{x+\tilde{s}} \int_x^{x+\tilde{s}} \* \rho_3(x, y_1, y_2) \* dy_1 \* dy_2 \right) \* dx.
\end{equation}
As before one can show that $\rho_3(x, y_1, y_2)=K[x,y_1, y_2]= O(\tilde{s}^4), \ $ and thus
$ E \* N_2(L) \leq (\tilde{s})^6 \*L=o(1). \ $

Theorem 1 is proven.  Theorem 2 immediately follows from Theorem 1. Indeed, the event $ \{ \eta(L) >s \} $ is exactly the event that there are
no nearest spacings smaller than $ s \* L^{-\frac{1}{3}} $ between the particles in $[0, L].$

{\bf Acknowledgment} {\it Research was supported in part by the NSF Grant DMS-0405864.}

 \bibliographystyle{}
 \bibliography{}
 \def\cmp{{\it Commun. Math. Phys.} }

\printindex
\end{document}